\documentclass[11pt]{article}
\usepackage{amssymb,amsmath,latexsym}

\newtheorem{example}[equation]{Example}
\newtheorem{theorem}[equation]{Theorem} 
\newtheorem{corollary}[equation]{Corollary}
\newtheorem{lemma}[equation]{Lemma}
\newtheorem{proposition}[equation]{Proposition}
\newtheorem{definition}[equation]{Definition}

\numberwithin{equation}{section}

\newcommand{\A}{\mathfrak{a}}
\newcommand{\fA}{\mathfrak{A}}
\newcommand{\I}{\mathcal{I}}
\newcommand{\J}{\mathcal{J}}
\newcommand{\qed}{\hfill $\Box$}
\newcommand{\Z}{\mathbb{Z}}
\newcommand{\proof}{{\bf Proof\ \ }}
\newcommand{\fg}{\mathfrak{g}}
\newcommand{\fgtilde}{\widetilde{\fg}}
\newcommand{\fghat}{\widehat{\fg}}
\newcommand{\ot}{\otimes}
\newcommand{\ol}{\underline}
\newcommand{\rhobr}{\}_\rho}

\newcommand{\Azh}{\widehat{\A}_0}
\newcommand{\gab}{\underline{\ga}}

\newcommand{\gbb}{{\underline{\gb}}}
\newcommand{\ggb}{\underline{\gamma}}
\newcommand{\wab}{w_{\underline{\ga}}}
\newcommand{\lbb}{{\lambda_{\underline{\beta}}}}
\newcommand{\wgbs}{w_{\underline{\gamma}}^*}
\newcommand{\Cq}{{\mathbb C}_q}
\newcommand{\ga}{\alpha}
\newcommand{\gb}{\beta}
\newcommand{\bc}{{\bf c}}
\newcommand{\ep}{\epsilon}
\newcommand{\gd}{\delta}
\newcommand{\gl}{\lambda}
\newcommand{\wa}{w_\alpha}

\newcommand{\lb}{\lambda_\beta}

\newcommand{\glA}{\mathfrak{gl}_N(A)}

\newcommand{\glK}{\mathfrak{gl}_N({\mathbb K})}

\newcommand{\Ares}{\A_{\rm res}}
\newcommand{\Vres}{V_{\rm res}}
\newcommand{\wgb}{w_{\underline{\gamma}}}
\newcommand{\wmbs}{w_{\underline{\mu}}^*}
\newcommand{\wbbs}{w_{\ol{\beta}}^*}

\title{Bosonic and Fermionic Representations \\
of Lie Algebra Central Extensions}

\date{May 28, 2004}

\author{Michael Lau\thanks{The author wishes to thank Professor G.~Benkart for her many helpful suggestions, and Professor B.N.~Allison for his comments on a previous draft.
This work was partially supported by N.S.A. grant MDA 904-03-1-0068.\newline
E-mail: lau@math.wisc.edu} \\ Department of Mathematics, University of Wisconsin \\ Madison, Wisconsin 53706, U.S.A.}

\begin{document}
\maketitle

\begin{small}
\noindent 
{\bf Abstract.}  Given any representation of an arbitrary Lie algebra $\fg$ over a field ${\mathbb K}$ of characteristic 0, we construct representations of a central extension of $\fg$ on bosonic and fermionic Fock space.  The method gives an explicit formula for a (sometimes trivial) 2-cocycle in $H^2(\fg ; {\mathbb K})$.  We illustrate these techniques with several concrete examples.

\bigskip



\end{small}

\vskip.25truein
\section{Introduction}

Clifford and Weyl algebras have natural representations on exterior and symmetric algebras, respectively.  In the early 1980s, I.B. Frenkel, V.G. Kac, and D.H. Peterson (\cite{Fr1},\cite{KP}) explicitly constructed the orthogonal affine Lie algebra $\widehat{\mathfrak{so}_N}$ using the quadratic elements of a Clifford algebra $C$.  These elements were viewed as quadratic operators on a certain highest weight $C$-module, called fermionic Fock space.  A.J. Feingold and I.B. Frenkel \cite{FF} later gave an analogous construction of the symplectic affine Lie algebra $\widehat{\mathfrak{sp}_{2N}}$ from quadratic elements of a Weyl algebra $A$ acting on bosonic Fock space, a highest weight $A$-module.  The natural inclusion of $\widehat{\mathfrak{gl}_N}$ into both $\widehat{\mathfrak{so}_{2N}}$ and $\widehat{\mathfrak{sp}_{2N}}$ meant that the affine Lie algebra of type $A_{N-1}^{(1)}$ admits a uniform construction of both fermionic and bosonic modules.  Notably, the bosonic construction (which has level $-1$) was the first known construction of a nonstandard irreducible highest weight representation\footnote{The {\it level} $k$ of a representation of an affine Lie algebra is the constant by which the canonical central element acts. If $k$ is anything other than a nonnegative integer, the representation is said to be {\it nonstandard}.} for an affine Lie algebra.  Seligman \cite{Sel} later modified the Feingold-Frenkel construction to produce a large class of irreducible integrable highest-weight modules for $\widehat{\mathfrak{so}_{N}}$ and $\widehat{\mathfrak{sp}_{2N}}$.  The fermionic construction (for types $A$ and $D$) is isomorphic to the vertex operator construction (in \cite{Fr2},\cite{tKvL} for instance), giving a boson-fermion correspondence.  A generalization of this boson-fermion correspondence appears in \cite{La}.

In 2002, Y. Gao \cite{G} used techniques similar to those of Feingold and Frenkel \cite{FF} to construct bosonic and fermionic representations of the extended affine Lie algebra $\widetilde{\mathfrak{gl}_N({\mathbb C}_q)}$, where $\Cq$ is the quantum torus in two variables.  This was accomplished by defining an interesting module for $\widehat{\mathfrak{gl}_N(\Cq)}$, a central extension of $\mathfrak{gl}_N(\Cq)$.

The Feingold-Frenkel construction for $\widehat{\mathfrak{gl}_N({\mathbb C})}$ and Gao's construction for $\widehat{\mathfrak{gl}_N(\Cq)}$ are special cases of a more general phenomenon addressed in this paper.  Both constructions define a Weyl or Clifford algebra ${\mathfrak A}$ from generators that we view as basis elements of modules for the Lie algebras.  These modules, ${\mathbb C}^N\ot{\mathbb C}[t,t^{-1}]$ and its dual module $({\mathbb C}^N\ot{\mathbb C}[t,t^{-1}])^*$, can be thought of as natural modules for the Lie algebras $\mathfrak{gl}_N({\mathbb C}[t,t^{-1}])$ and ${\mathfrak{gl}}_N(\Cq)$.  Selecting some of these generators to be ``positive'', Feingold-Frenkel and Gao define an associative subalgebra ${\mathfrak A}^+$ of ${\mathfrak A}$.  Their Fock spaces are trivial ${\mathfrak A}^+$-modules induced to all of ${\mathfrak A}$.
        
Our construction replaces ${\mathfrak{gl}_N({\mathbb C}[t,t^{-1}])}$ and ${\mathfrak{gl}_N(\Cq)}$ with {\em any} Lie algebra $\fg$ over a field ${\mathbb K}$ of characteristic $0$, and replaces ${\mathbb C}^N\ot{\mathbb C}[t,t^{-1}]$ with an {\em arbitrary} $\fg$-module $W$.  We use Weyl- (resp. Clifford-) type relations to generate a unital associative algebra $\A$.  This algebra is constructed from a basis of $W\oplus W'$ where $W'$ is the $\fg$-submodule of the dual space $W^*$ generated by the restricted dual $\hbox{\rm Span}_{\mathbb K}\{ \wa^*\ |\ \ga\in\I\}$, where $\{\wa\ |\ \ga\in\I\}$ is a basis for $W$.  The generators for $\A$ are partitioned into ``positive'' and ``negative'' elements, and the division is used to define a vacuum vector and an induced module for $\A$, called bosonic (resp. fermionic) Fock space $V$.  Some care is needed so that $V$ remains a well-defined module under the action of quadratic operators $f_x$ defined for each element $x$ in $\fg$.  We treat the operators $f_x$ as elements of a completion of $\A$.  The fact that the completion is itself an associative algebra simplifies some of the most difficult computations of Feingold-Frenkel and Gao.

%Namely, {\it any} representation of {\it any} Lie algebra can be used to uniformly construct bosonic and fermionic representations of a central extension of that Lie algebra.

The assignment $x\mapsto f_x$ extends to a representation of a certain one-dimensional central extension $\fgtilde$ of $\fg$.  For Feingold-Frenkel, $\fgtilde$ is simply the affine Lie algebra $\widehat{\mathfrak{gl}_N({\mathbb C})}$ and for Gao, $\fgtilde$ is a homomorphic image of $\widehat{\mathfrak{gl}_N(\Cq)}$.  The representation he constructs for 
$\widehat{\mathfrak{gl}_N(\Cq)}$ is the pullback of the Fock representation we construct for $\fgtilde$.

When this representation is faithful, it gives an embedding of $\fgtilde$ into a Lie algebra of ``infinite matrices'', analogous to that found by Kac-Peterson \cite{KP} for $\widehat{\mathfrak{gl}_N({\mathbb C})}$.  It also affords a simple formula for the 2-cocycle defining the central extension, and we use it to explicitly compute the 2-cocycle of the (universal) central extension of the Lie algebra of $N\times N$ matrices over the ring of differential operators of the form $\sum_{n=0}^\infty f_n(t)(\frac{d}{dt})^n$ (where $f_n(t)\in {\mathbb C}[t,t^{-1}]$ is $0$ for $n\gg 0$).

\section{Bosonic and Fermionic Realizations}

Assume $\rho=1$ or $-1$.  If $\rho=-1$ (resp. $+1$), we call the resulting constructions {\it bosonic} (resp. {\it fermionic}).  For elements $a,b$ of any associative algebra $A$, let $\{ a,b\rhobr=ab+\rho ba$.  Note that $\{a,b\rhobr=\rho\{ b,a\rhobr$, and $[ab,c]=a\{ b,c\rhobr -\rho\{ a,c\rhobr b$ for $a,b,c\in A$, where $[a,b]$ is the usual commutator $ab-ba$.

\begin{definition}
Let $\fg$ be a Lie algebra over a field ${\mathbb K}$ of characteristic 0, and let $W$ be an arbitrary $\fg$-module with ${\mathbb K}$-basis $\mathfrak{B}=\{ w_\ga \ |\ \ga\in\I\}$.  The universal enveloping algebra $\mathcal{U}(\fg )$ has a natural action on the dual module $W^*=\hbox{\rm Hom}_{\mathbb K}(W,{\mathbb K})$ coming from the $\fg$-action $(x.\gl )(w)=-\gl (x.w)$ for all $x\in\fg,\ \gl\in W^*,\ w\in W$.  Let $W'=\sum_{\ga\in\I}\mathcal{U}(\fg ).w_\ga^*$ where the linear functionals $\wa^*:\ W\rightarrow {\mathbb K}$ are defined by $\wa^*(w_\gb)=\gd_{\ga,\gb}\ \forall \ga,\gb\in\I$.  Note that if $W$ is infinite-dimensional, it is possible that $W'\supsetneq \hbox{\rm Span}_{\mathbb K}\{ \wa^*\ |\ \ga\in\I\}$, so we fix a ${\mathbb K}$-basis $\mathfrak{B}'=\{\wa^*,\lb\ |\ \ga\in\I,\ \gb\in\I '\}$ for $W'$.  The choice of $\mathfrak{B}$, $\mathfrak{B}'$, and a subset $\J\subseteq\I$ is called a {\bf realization} of $(\fg,W)$ if $\hbox{\rm card}\{\ga\in\J\ |\ x.\wa\notin S_\J\}<\infty$ for each $x\in\fg$, where $S_\J=\hbox{\rm Span}_{\mathbb K}\{w_\gb\ |\ \gb\in\J\}$.
\end{definition}

Realizations always exist---for example, any finite subset $\J\subseteq\I$ will trivially satisfy the ``finiteness condition'', $\hbox{\rm card}\{\ga\in\J\ |\ x.\wa\notin S_\J\}<\infty$.  The purpose of a realization is to define an associative algebra $\A$, with ``positive'' and ``negative'' parts, $\A^+$ and $\A^-$, for later use in constructing bosonic (resp. fermionic) Fock space.  

Given a realization of a Lie algebra $\fg$ and representation $W$, let $\A=\A (\fg,W,\I,\I ',\rho )$ be the (unital) associative algebra generated by 

\noindent
$\{ w_\ga ,  w_\ga^*,\gl_\gb\ | \ \ga\in\I,\gb\in\I '\},$ modulo the relations

\begin{itemize}
\item[{\rm (R1)}] $\{v,w\rhobr=\{\gl,\eta\rhobr=0$
\item[{\rm (R2)}] $\{\gl,w\rhobr=\gl(w)$
\end{itemize}

\noindent
for all $v,w\in W$, $\gl,\eta\in W'$.  Let $\A^+\subseteq\A$ denote the (unital) subalgebra generated by those elements $\wa\in \mathfrak{B}$, $\gl\in \mathfrak{B}'$ such that $\ga\in\J$ and $\gl (w_\gb)=0$ for all $\gb\in\J$.  Likewise, let $\A^-$ be the (unital) subalgebra generated by $\{ \wa\in\mathfrak{B},\gl\in\mathfrak{B}'\ |\ \ga\in\I\setminus\J,\mbox{ and }\gl (w_\gb )\neq 0\mbox{ for some }\gb\in\J\}$.  We will sometimes use the {\it restricted algebra} $\Ares$, defined as the (unital) subalgebra generated by $\{\wa,\wa^*\ |\ \ga\in\I\}$.  Analogously, $\Ares^+$ and $\Ares^-$ are the (unital) subalgebras generated by $\{\wa\in\mathfrak{B},w_\gb^*\in\mathfrak{B}'\ |\ \ga\in\J,\ \gb\in\I\setminus\J\}$ and $\{\wa\in\mathfrak{B}, w_\gb^*\in\mathfrak{B}'\ |\ \ga\in\I\setminus\J,\ \gb\in\J\}$, respectively.  Note that in many interesting cases (see Examples \ref{ex6.1}-\ref{ex6.7} below), $\I'=\emptyset$, so $\Ares=\A$ and $\Ares^\pm=\A^\pm$.\footnote{An example where $\I'\neq\emptyset$ is where $\fg={\rm Span}_{\mathbb K}\{\sum_{n=0}^\infty f_n(t)(\frac{d}{dt})^n\ |\ f_n(t)\in {\mathbb K}[t]\}$, viewed as a Lie algebra of infinite series of differential operators on $W={\mathbb K}[t]$.}

We will use the multiindex notation $\wab=w_{\ga_1}w_{\ga_2}\cdots w_{\ga_r},$ 

\noindent
$\wab^*=w_{\ga_1}^*w_{\ga_2}^*\cdots w_{\ga_r}^*$, and $\lbb=\gl_{\gb_1}\gl_{\gb_2}\cdots\gl_{\gb_s}$
where $\ol{\ga}=(\ga_1,\ldots ,\ga_r)$ and $\ol{\gb}=(\gb_1,\ldots ,\gb_s)$.  Fix a total ordering $\preceq$ on $\I$ and $\I'$.  We say that $\wab\gl_{\gbb}\wgbs\neq 0$ is a {\it standard monomial} if $\ga_i\preceq\ga_j$, $\gb_i\preceq\gb_j$, and $\gamma_i\preceq\gamma_j$whenever $i<j$.  By relations (R1) and (R2), the standard monomials span $\A$.  The {\it length} $\ell(\ol{\ga})$ of a multiindex $\ol{\ga}=(\ga_1,\ldots, \ga_r)$ is $r$, the number of its entries.  If $\ol{\ga}=\emptyset$, we say that $\ell(\gab)=0$, and we let $\wab=\gl_{\ol{\ga}}=\wab^*=1$. 

Let $\widehat{\A}$ be the ${\mathbb K}$-vector space of linear combinations of (possibly infinitely many) distinct standard monomials $\wab\lbb\wgbs$ of $\A$.  That is, the elements of $\widehat{\A}$ are those that can be expressed as (possibly infinite) sums $\sum_{\gab,\gbb,\ggb}c_{\gab\gbb\ggb}\wab\lbb\wgbs$ where $c_{\gab\gbb\ggb}\in {\mathbb K}$, and the sum runs over all finite multiindices $\gab, \ggb$ with entries in $\I$ and $\gbb$ with entries in $\I'$ such that $\wab\lbb\wgbs$ is a standard monomial.  Elements written in this way (as linear combinations of distinct standard monomials) are said to be in a {\it standard form}.  For the remainder of the paper, we will restrict our attention to the subspace $\Azh$ consisting of those elements $\sum_{\gab,\gbb,\ggb}c_{\gab\gbb\ggb}\wab\lbb\wgbs\in\widehat{\A}$ with the following property: 

\begin{itemize}
\item[{\rm (P1)}] For each $\ggb$, there are only finitely many pairs $(\gab,\gbb)$ such that $c_{\gab\gbb\ggb}\neq 0$.
\end{itemize}

{\lemma Term-by-term multiplication gives $\Azh$ the structure of an associative algebra.}

\vspace{0.1in}
\noindent
\proof  It is enough to show that this multiplication is well-defined and the resulting products are in $\Azh$.  Associativity then follows immediately from the associativity of $\A$.  Fix a multiindex $\ol{\zeta}=(\zeta_1,\ldots ,\zeta_t)$ of elements from $\I$.  Let $T=\displaystyle{\sum_{\gab,\gbb,\ggb} c_{\gab\gbb\ggb}\wab\lbb\wgbs}$ and $T'=\displaystyle{\sum_{\ol{\ga'},\ol{\gb'},\ol{\gamma'}} d_{\ol{\ga'}\ol{\gb'}\ol{\gamma'}}w_{\ol{\ga'}}\gl_{\ol{\gb'}}w_{\ol{\gamma'}}^*  }$ be elements of $\Azh$ written in a standard form.

We consider the products $P=P_{\gab\gbb\ggb\ol{\ga'}\ol{\gb'}\ol{\gamma'}}=c_{\gab\gbb\ggb}\wab\lbb\wgbs d_{\ol{\ga'}\ol{\gb'}\ol{\gamma'}}w_{\ol{\ga'}}\gl_{\ol{\gb'}}w_{\ol{\gamma'}}^*$ which, when written in a standard form, contain a nonzero multiple of a standard monomial of the form $w_{\ol{\eta}}\gl_{\ol{\mu}}w_{\ol{\zeta}}^*$.  Note that $w_{\gamma_i'}^*\in\{ w_{\zeta_j}^*\ |\ 1\leq j\leq t\}$ for all $w_{\gamma_i'}^*$ occurring in the expression for $P$.  Thus, there are only finitely many possible $w_{\ol{\gamma'}}^*$ that may occur.  For each such $\ol{\gamma'}$, there are only finitely many $\ol{\ga'}=(\ga_1',\ldots ,\ga_r')$, $\ol{\gb'}=(\gb_1',\ldots,\gb_s')$ for which $d_{\ol{\ga'}\ol{\gb'}\ol{\gamma'}}\neq 0$, since $T'\in\Azh$.  For each of these (finitely many) possible triples $(\ol{\ga'},\ol{\gb'},\ol{\gamma'})$ occurring in the expression for $P$, we see that every $w_{\gamma_i}^*\in\{ w_{\ga_j'}^*, w_{\zeta_k}^*\ |\ 1\leq j\leq r,\ 1\leq k\leq t\}$.  (Any $w_{\gamma_i}^*\notin\{ w_{\ga_j'}^*\ |\ 1\leq j\leq r\}$ would commute with $w_{\ol{\ga'}}\gl_{\ol{\gb'}}$ and thus contribute a $w_{\gamma_i}^*$ term to the expression $w_{\ol{\zeta}}^*$, which is impossible unless $w_{\gamma_i}^*\in\{ w_{\zeta_k}^*\ |\ 1\leq k\leq t\}$.)  For each of these finitely many possible $\ol{\gamma}$, there are only finitely many $(\ol{\ga},\ol{\gb})$ such that $c_{\gab\gbb\ggb}\neq 0$, so there are only finitely many $(\gab,\gbb,\ggb,\ol{\ga'},\ol{\gb'},\ol{\gamma'})$ such that a standard form expression for $P$ contains a term of the form $kw_{\ol{\eta}}\gl_{\ol{\mu}}w_{\ol{\zeta}}^*$, where $k$ is a nonzero scalar.  That is, the product $TT'$ is well-defined, and in $\Azh$.\qed

\vspace{0.15in}
For each $x\in\fg$ and $\ga\in \I$, write $x.\wa=\displaystyle{\sum_{\gamma\in\I}x_\gamma^\ga w_\gamma}$ with $x_\gamma^\ga\in {\mathbb K}$.  Note that $x_\gamma^\ga=0$ for all but finitely many $\gamma$.  Define the normal ordering $:\sum c_{\ga\gb}\wa w_\beta^*:=\sum c_{\ga\gb}:\wa w_\beta^*:$ where $c_{\ga\gb}\in {\mathbb K}$ and

\begin{equation}\label{eqn1}
:\wa w_\beta^*:=\left\{\begin{array}{ll}
-\rho w_\beta^*\wa & \mbox{if }\ga=\gb\in\J \\
\wa w_\beta^* & \mbox{otherwise.}
\end{array}\right.
\end{equation}
Let 
\begin{equation}\label{eqn2}
f_x=\sum_{\ga\in\I}:(x.\wa)\wa^*:.
\end{equation}
  Note that $f_x=\sum_{\ga,\gamma\in\I}x_\gamma^\ga:w_\gamma\wa^*:$ is a well-defined member of $\Azh$ for each $x\in\fg$.  Thus, we may conduct the multiplications of the following two lemmas within $\Azh$.

{\lemma\label{lm2.2} For every $x\in\fg$, $\eta\in\I$, and $\gl\in W'$, 

\begin{itemize}
\item[{\rm (i)}]
$[f_x,w_\eta]=\sum_{\gamma\in\I}x_\gamma^\eta w_\gamma=x.w_\eta$

\item[{\rm (ii)}]
$[f_x,\gl]=-\sum_{\ga,\gamma\in\I}x_\gamma^\ga\gl(w_\gamma)\wa^*=x.\gl$.
\end{itemize}}

\vspace{0.1in}
\noindent
\proof  If $\gamma = \ga\in\J$, then $:w_\gamma\wa^*:\ = -\rho\wa^* w_\gamma = w_\gamma\wa^*+1$.  Otherwise, $:w_\gamma\wa^*:\ =w_\gamma\wa^*$.  Therefore $\hbox{{\rm ad}}:w_\gamma\wa^*:\ =\hbox{{\rm ad }}w_\gamma\wa^*$, so

$$\begin{array}{ll}
[f_x,w_\eta]&=\sum_{\ga,\gamma\in\I}x_\gamma^\ga[:w_\gamma\wa^*:,w_\eta] \\
&=\sum_{\ga,\gamma\in\I}x_\gamma^\ga[w_\gamma\wa^*,w_\eta]\\
&=\sum_{\ga,\gamma\in\I}x_\gamma^\ga w_\gamma\{\wa^*,w_\eta\rhobr\\
&=\sum_{\gamma\in\I}x_\gamma^\eta w_\gamma=x.w_\eta.\\
\end{array}$$  

\noindent
Similarly, $$\begin{array}{ll}
[f_x,\gl]&=-\rho\sum_{\ga,\gamma\in\I}x_\gamma^\ga \{w_\gamma,\gl\rhobr\wa^*\\
&=-\sum_{\ga,\gamma\in\I}x_\gamma^\ga\gl(w_\gamma)\wa^*,\\
\end{array}$$ 

\noindent
and $$\begin{array}{ll}
(x.\gl)(\wa)&=-\gl(x.\wa)\\
&=-\sum_{\gamma\in\I}x_\gamma^\ga\gl(w_\gamma).\\
\end{array}$$  

\noindent
Thus $x.\gl=-\sum_{\ga,\gamma\in\I}x_\gamma^\ga\gl(w_\gamma)\wa^*=[f_x,\gl]$.\qed

{\lemma\label{lm2.3} $[f_x,f_y]=f_{[x,y]}+\Omega_{x,y}^\rho$ for some $\Omega_{x,y}^\rho\in\Azh$ such that $[\Omega_{x,y}^\rho,\A]=0$.}

\vspace{0.1in}

\noindent
\proof  By the Jacobi identity (on $\Azh$, viewed as a Lie algebra), if $u\in W$ or $W'$, then $[[f_x,f_y],u]=[f_x,[f_y,u]]-[f_y,[f_x,u]]=x.(y.u)-y.(x.u)=[x,y].u=[f_{[x,y]},u]$.  Thus $[f_x,f_y]-f_{[x,y]}$ commutes with the generators of $\A$, and hence centralizes $\A$.\qed

\section{Fock Space}

Define an augmentation map $\ep:\ \A^+\rightarrow {\mathbb K}$ by the rule that
$$\ep(\wab\lbb\wgbs)=\left\{\begin{array}{ll}
1 & \mbox{if }\ol{\ga}=\ol{\gb}=\ol{\gamma}=\emptyset \\
0 & \mbox{otherwise.}
\end{array}
\right.
$$
Let ${\mathbb K}v_0$ be the one-dimensional left $\A^+$-module given by $$\mu.v_0=\ep(\mu)v_0\mbox{ for }\mu\in\A^+.$$  
\noindent
Inducing to $\A$ gives the {\it Fock space} $V=V(\fg,W,\J,\I,\I',\rho)=\A\ot_{\A^+}{\mathbb K}v_0$.  

Analogously, we define the left $\Ares$-module
$\Vres=V(\fg,W,\J,\I,\rho)=\Ares\otimes_{\Ares^+}{\mathbb K}v_0$,
where $\Ares^+$ acts on ${\mathbb K}v_0$ by restriction of the $\A^+$-action.  We write $av_0$ to denote $a\otimes v_0$ in either of $V$ or $\Vres$.  Which module $av_0$ inhabits should be clear from context.

By the relations of $\A$, each $v\in\ V$ (resp. $\Vres$) can be written in the form $av_0$ for some $a\in\A^-$ (resp. $\Ares^-$).  From the relations of $\A$, we have the following useful formulas:

{\lemma\label{comp1} Let $\wa,w_\gamma,\wa^*,w_\gamma^*\in\Ares$.  Then
\begin{itemize}

\item[{\rm (i)}] $\wa w_\gamma=\left\{\begin{array}{ll}
w_\gamma \wa & \mbox{if }\rho=-1 \\
-w_\gamma \wa & \mbox{if }\rho=1; \\
\end{array}
\right.
$

\item[{\rm (ii)}] $\wa w_\gamma^*=\left\{\begin{array}{ll}
w_\gamma^* \wa-\gd_{\ga,\gamma}1 & \mbox{if }\rho=-1 \\
-w_\gamma^* \wa+\gd_{\ga,\gamma}1 & \mbox{if }\rho=1; \\
\end{array}
\right.
$

\item[{\rm (iii)}] $\wa^* w_\gamma=\left\{\begin{array}{ll}
w_\gamma \wa^*+\gd_{\ga,\gamma}1 & \mbox{if }\rho=-1 \\
-w_\gamma \wa^*+\gd_{\ga,\gamma}1 & \mbox{if }\rho=1; \\
\end{array}
\right.
$

\item[{\rm (iv)}] $\wa^* w_\gamma^*=\left\{\begin{array}{ll}
w_\gamma^* \wa^* & \mbox{if }\rho=-1 \\
-w_\gamma^* \wa^* & \mbox{if }\rho=1. \\
\end{array}
\right.
$

\end{itemize}
}\qed

\vspace{0.15in}
\noindent
Iterating Lemma \ref{comp1} gives

{\lemma\label{comp2}
Suppose $\wa,\wa^*\in\Ares$, and let 
$$0\neq\wgb\wmbs=w_{\gamma_1}^{m_1}w_{\gamma_2}^{m_2}\cdots w_{\gamma_r}^{m_r}w_{\mu_1}^{*n_1}w_{\mu_2}^{*n_2}\cdots w_{\mu_s}^{*n_s}\in\Ares$$ 

\noindent
where
each $m_i,n_j$ is a positive integer,\footnote{Note that if $\rho=1$, then each $m_i$ and $n_j$ is equal to $1$.} the $\gamma_i$ are all pairwise distinct, and the $\mu_j$ are all pairwise distinct.

\noindent
Then

\noindent
$\hbox{{\rm (i)\ \  }} \displaystyle{\wa\wgb\wmbs}\ =\ \displaystyle{(-\rho)^{\ell(\ggb)+\ell(\ol{\mu})}\wgb\wmbs\wa}$

$$\hspace{0.15in}-\ \displaystyle{\wgb\sum_{j=1}^s(-\rho)^{\ell(\ggb)+j+1}n_j\gd_{\ga,\mu_j}w_{\mu_1}^{*n_1}\cdots w_{\mu_{j-1}}^{*n_{j-1}}w_{\mu_j}^{*n_j-1}w_{\mu_{j+1}}^{*n_{j+1}}\cdots w_{\mu_s}^{*n_s}}$$

\noindent
$
\hbox{{\rm (ii)\ \  }} \displaystyle{\wa^*\wgb\wmbs}\ =\ \displaystyle{(-\rho)^{\ell(\ggb)+\ell(\ol{\mu})}\wgb\wmbs\wa^*}$
$$\hspace{0.15in}+\ \displaystyle{\sum_{i=1}^r(-\rho)^im_i \gd_{\ga,\gamma_i}w_{\gamma_1}^{m_1}\cdots w_{\gamma_{i-1}}^{m_{i-1}}w_{\gamma_i}^{m_i-1}w_{\gamma_{i+1}}^{m_{i+1}}\cdots w_{\gamma_r}^{m_r}\wmbs}.$$
}\qed

\vspace{0.15in}
Thus if $\wgb\wmbs\in\A^-$, the action of $\wa$ (resp. $w_\gb^*$)$\in\A^+$ on $\wgb\wmbs v_0$ is, up to a factor of $\pm 1$, $\frac{\partial}{\partial \wa^*}$ (resp. $\frac{\partial}{\partial w_\gb}$):

{\lemma\label{comp3} Let $\wa,w_\gb^*\in\A^+$, and suppose $\wgb\wmbs$ is as in Lemma \ref{comp2}.  Then

\noindent
$
\hbox{{\rm (i)\ \  }} \displaystyle{\wa\wgb\wmbs}v_0$

$$=\ -\displaystyle{\sum_{j=1}^s(-\rho)^{\ell(\ggb)+j+1}n_j \gd_{\ga,\mu_j}\wgb w_{\mu_1}^{*n_1}\cdots w_{\mu_{j-1}}^{*n_{j-1}}w_{\mu_j}^{*n_j-1}w_{\mu_{j+1}}^{*n_{j+1}}\cdots w_{\mu_s}^{*n_s}}v_0$$

\noindent
$
\hbox{{\rm (ii)\ \  }} \displaystyle{\wa^*\wgb\wmbs}v_0=\displaystyle{\sum_{i=1}^r(-\rho)^im_i \gd_{\ga,\gamma_i}w_{\gamma_1}^{m_1}\cdots w_{\gamma_{i-1}}^{m_{i-1}}w_{\gamma_i}^{m_i-1}w_{\gamma_{i+1}}^{m_{i+1}}\cdots w_{\gamma_r}^{m_r}\wmbs}v_0.$}

\qed

{\corollary\label{monomaker} Suppose $\ell(\ol{\ga})+\ell(\ol{\gb})\geq\ell(\ol{\gamma})+\ell(\ol{\mu})$, $w_{\ol{\ga}}\wbbs\in\A^*$, and $\wgb\wmbs$ is as in Lemma \ref{comp2}.  Then 
$$\displaystyle{w_{\ol{\ga}}w_{\ol{\gb}}^*\wgb\wmbs v_0=k(\ggb,\ol{\mu})\gd_{\ol{\ga},\ol{\mu}}\gd_{\ol{\gb},\ol{\gamma}}m_1!m_2!\cdots m_r!n_1!n_2!\cdots n_s!v_0}$$

\noindent
where $k(\ggb,\ol{\mu})=\pm 1$.\qed
}

{\proposition\label{V_res} $\Vres$ is a simple $\Ares$-module with ${\mathbb K}$-basis 
$$\mathcal{B}=\{\wab\wbbs v_0\ |\ \wab\wbbs\in\A^-\ \hbox{ is a standard monomial}\}.$$

}

\vspace{0.1in}
\noindent
\proof By the observations before Lemma \ref{comp1}, the elements of $\mathcal{B}$ span $\Vres$.  Let $v\in\Vres$.  Write $v=\sum_{\ol{\ga},\ol{\gb}}c_{\ol{\ga}\ol{\gb}}\wab\wbbs v_0=0$ where each $\wab\wbbs$ is a standard monomial in $\A^-$.  Choose multiindices $\ol{\ga}'$ and $\ol{\gb}'$ so that $w_{\ol{\ga}'}w_{\ol{\gb}'}\in\A^-$ and $\ell(\ol{\ga}')+\ell(\ol{\gb}')=\hbox{{\rm max}}\{\ell(\ol{\ga})+\ell(\ol{\gb})\ |\ c_{\ol{\ga}\ol{\gb}}\neq 0\}$.  Then by Corollary \ref{monomaker}, 

\begin{equation}\label{e1}
w_{\ol{\gb}'}w_{\ol{\ga}'}^*v=kc_{\ol{\ga}'\ol{\gb}'}v_0
\end{equation}

\noindent
for some $k\in {\mathbb K}\setminus 0$.  If $v\neq 0$, then $c_{\ol{\ga}'\ol{\gb}'}$ can be taken to be nonzero.  Hence every nonzero submodule $V'$ contains $v_0$.  But since $v_0$ generates $\Vres$, $V'$ must be $\Vres$, so $\Vres$ is simple.  

If $v=0$, then (\ref{e1}) reads $0=w_{\ol{\gb}'}w_{\ol{\ga}'}^*v=kc_{\ol{\ga}'\ol{\gb}'}v_0$.  Hence $c_{\ol{\ga}'\ol{\gb}'}=0$, and by the maximality of $\ell(\ol{\ga}')+\ell(\ol{\gb}')$, we see that every $c_{\ol{\ga}\ol{\gb}}$ is zero.  Therefore the elements of $\mathcal{B}$ are linearly independent and form a ${\mathbb K}$-basis for $\Vres$.\qed

{\proposition\label{Kvzerocrit} Let $v\in\Vres$.  Then $v\in {\mathbb K}v_0$ if and only if 

$\big((1-\ep)(\Ares^+)\big).v=0.$}

\vspace{0.1in}
\noindent
\proof If $v\in {\mathbb K}v_0$, then $\big((1-\ep)(\Ares^+)\big)v=0$, by the definition of $\Vres$.  Conversely, assume that $\big((1-\ep)(\Ares^+)\big)v=0$.  Write $v=\sum_{\ol{\ga},\ol{\gb}}c_{\ol{\ga}\ol{\gb}}\wab\wbbs v_0=0$ where each $c_{\ol{\ga}\ol{\gb}}\in {\mathbb K}$, and $\wab\wbbs$ is an element of $\Ares^-$ in standard form.  Choose multiindices $\ol{\ga}',\ol{\gb}'$ such that $\ell(\ol{\ga}')+\ell(\ol{\gb}')=\hbox{{\rm max}}\{\ell(\ol{\ga})+\ell(\ol{\gb})\ |\ c_{\ol{\ga}\ol{\gb}}\neq 0\}$ and $c_{\ol{\ga}'\ol{\gb}'}\neq 0$.  By Corollary \ref{monomaker}, $w_{\ol{\gb}'}w_{\ol{\ga}'}^*v=kc_{\ol{\ga}'\ol{\gb}'}v_0$ for some $k\in {\mathbb K}\setminus 0$.  But $w_{\ol{\ga}'}w_{\ol{\gb}'}^*\in\Ares^-$, so by the definition of $\Ares^+$,  $w_{\ol{\gb}'}w_{\ol{\ga}'}^*\in\Ares^+$.  If $\ol{\ga}'$ and $\ol{\gb}'$ are not both equal to $\emptyset$, then $w_{\ol{\gb}'}w_{\ol{\ga}'}^*v=\big((1-\ep)(w_{\ol{\gb}'}w_{\ol{\ga}'}^*)\big)v=0$, a contradiction.  Hence $\ol{\ga}'=\ol{\gb}'=\emptyset$, so we are done by the maximality of $\ell(\ol{\ga}')+\ell(\ol{\gb}')$.\qed

\vspace{0.15in}

Modules like $V$ play an important role in statistical mechanics, where they represent the ``space of states'' for a given system.  The cyclic vector $v_0$ is viewed as a vacuum, and the element 
$$w_{\ga_1}^{m_1}w_{\ga_2}^{m_2}\cdots w_{\ga_r}^{m_r}\gl_{\gb_1}^{n_1}\cdots \gl_{\gb_s}^{n_s} w_{\gamma_1}^{*{q_1}} \cdots w_{\gamma_t}^{*{q_t}}v_0$$
(with every $m_i,n_j,q_k$ a nonnegative integer, and $w_{\ga_i},\gl_{\gb_j},w_{\gamma_k}^*\in\A^-$) corresponds to a state with $m_i$ particles in state $w_{\ga_i}$, $n_j$ particles in state $\gl_{\gb_j}$, and $q_k$ particles in state $w_{\gamma_k}^*$.  When $\rho=-1$ (resp. $+1$), $V$ is called {\it bosonic} (resp. {\it fermionic}) {\it Fock space}, since the ``particles'' in $V$ satisfy Bose-Einstein (resp. Fermi-Dirac) occupancy statistics.

Due to the normal ordering and the ``finiteness condition'' in our definition of a realization, the elements $f_x$ defined in \S 2 have a well-defined left-multiplication action on $V$.  This follows easily from Lemma \ref{lm2.2} and the fact that all but finitely many monomials in any $f_x$ act as $0$ on the vacuum vector $v_0$.  We can therefore interpret the elements $f_x$ as operators on the Fock space $V$.

Moreover, for any $\ga,\gb,\gamma,\eta\in\I$, the bracket $[:w_\ga w_\gb^*:\, , :w_{\gamma}w_\eta^*:]$ is an element of $\hbox{Span}_{\mathbb K}\{:w_\mu w_\nu^*:\ |\ \mu,\nu\in\I\}\oplus{\mathbb K}1$, so $[f_x,f_y]$ may be written in the form $\sum_{\mu,\nu\in\I}c_{\mu\nu}:w_\mu w_\nu^*:$ for some $c_{\mu\nu}\in{\mathbb K}$.  Since $f_x$ and $f_y$ are given by a realization $(\fg,W,\J,\I,\I')$, the ``finiteness condition'' ensures that there are only finitely many $(\mu,\nu)\in(\I\setminus\J)\times\J$ such that $c_{\mu\nu}\neq 0$.  Thus $[f_x,f_y]v_0\in\Vres$.  Now by Lemma \ref{lm2.3}, $\Omega_{x,y}^\rho v_0\in\Vres$.  But since $\Omega_{x,y}^\rho$ commutes with the elements of $\A$, we see that
$$\big((1-\ep)(\Ares^+)\big)\Omega_{x,y}^\rho v_0 =\Omega_{x,y}^\rho\big((1-\ep)(\Ares^+)\big)v_0=0$$
Hence by Proposition \ref{Kvzerocrit}, $\Omega_{x,y}^\rho v_0\in {\mathbb K}v_0$.  Writing $\Omega_{x,y}^\rho v_0=c_{x,y}^\rho v_0$ with $c_{x,y}^\rho\in {\mathbb K}$ gives

{\corollary\label{LA} $\hbox{{\rm Span}}_{\mathbb K}\{ f_x\ |\ x\in\fg\}\oplus {\mathbb K}e$ is a Lie algebra with bracket defined by $[f_x, f_y]=f_{[x,y]}+c_{x,y}^\rho e$ and $[e,f_x]=0$ for all $x,y\in\fg$.}\qed

\section{Central Extensions}

A {\it central extension} of $\fg $ is a Lie algebra $\fgtilde$ and an epimorphism $\phi:\ \fgtilde\rightarrow\fg$ with $\ker\phi$ contained in the center $Z(\fgtilde)$ of $\fgtilde$.  Given two central extensions $(\fghat,\pi)$ and $(\fgtilde,\phi)$ of $\fg$, a {\it morphism} (from $(\fghat,\pi)$ to $(\fgtilde,\phi)$) in the category of central extensions is a Lie algebra homomorphism $\mu:\ \fghat\rightarrow\fgtilde$ such that $\phi\mu=\pi$.  We say that $(\fghat,\pi)$ and $(\fgtilde,\phi)$ are {\it isomorphic} if the morphism $\mu$ is a bijection.  The central extension $(\fghat,\pi)$ is {\it universal} if there is a unique morphism from it to every other central extension of $\fg$.

Isomorphism classes of one-dimensional central extensions are in bijective correspondence with cohomology classes in $H^2(\fg;{\mathbb K})$.  In particular, each class $[c]\in H^2(\fg;{\mathbb K})$ determines a central extension $\fgtilde = \fg\oplus {\mathbb K}\bc$ with the bracket $[\cdot,\cdot]\widetilde{\ }$ in $\fgtilde$ given by $[x,y]\widetilde{\ }=[x,y]+c(x,y)\bc$ for $x,y\in\fg$, where $\bc$ is central and $c$ is a representative of the class $[c]$.  Conversely, if $c:\ \fg\times\fg\rightarrow {\mathbb K}$ is ${\mathbb K}$-bilinear and $\fgtilde=\fg\oplus {\mathbb K}\bc$ is a Lie algebra under the bracket $[\cdot,\cdot]\widetilde{\ }$ defined above, then $c$ is a representative of a cohomology class $[c]$ in $H^2(\fg;{\mathbb K})$.  Direct computation from the chain complex defining $H^2(\fg;{\mathbb K})$ shows that $[c]\in H^2(\fg;{\mathbb K})$ if and only if $c$ satisfies the {\it 2-cocycle identities:}

(i) $c(x,y)=-c(y,x)$ and

(ii) $c\big(x,[y,z]\big)+c\big(y,[z,x]\big)+c\big(z,[x,y]\big)=0$.

{\theorem\label{thm4.1} For every $x,y\in\fg$, let $c(x,y)=\rho c_{x,y}^\rho$.  Then $[c]\in H^2(\fg;{\mathbb K})$.}

\vspace{0.1in}
\noindent
\proof  Since $x\mapsto f_x$ is a linear map, the space $S=\hbox{{\rm Span}}_{\mathbb K}\{ f_x\ |\ x\in\fg\}$ is a Lie algebra under the bracket $[f_x,f_y]_S=f_{[x,y]}$.  Then by Corollary \ref{LA}, the map $(f_x,f_y)\mapsto c_{x,y}^\rho$ is a 2-cocycle for $S$.  Finally, $[c]\in H^2(\fg;{\mathbb K})$ since $x\mapsto f_x$ is a homomorphism from $\fg$ to $S$, and $c(x,y)$ is a constant multiple of $c_{x,y}^\rho$.\qed

{\corollary\label{cor4.1} Let $\fgtilde=\fg\oplus {\mathbb K}\bc$ with $[x,y]\widetilde{\ }=[x,y]+c(x,y)\bc$ where $\bc$ is central, and let $\pi:\ \fgtilde\twoheadrightarrow\fg$ be the canonical projection.  Then $(\fgtilde,\pi)$ is a central extension of $\fg$.}\qed

{\theorem\label{thm4.2} The action $\bc.v=\rho v$, $x.v=f_xv$ for $x\in\fg$, $v\in V$ gives a representation of $\fgtilde $ on $V$.}

\vspace{0.1in}
\noindent
\proof  The map $x\mapsto f_x$, $\bc\mapsto \rho 1$ extends to a linear transformation on $\fgtilde$, so it suffices to check that 

\begin{eqnarray*}
[x,y]\widetilde{\ }.v&=&\big([x,y]+c(x,y)\bc\big).v\\ 
&=&f_{[x,y]}v+\rho c(x,y)v\\
&=&\big([f_x,f_y]-\Omega_{x,y}^\rho\big)v+\rho c(x,y) v\\
&=&[f_x,f_y]v\\
&=&f_xf_yv-f_yf_xv\\
&=&x.(y.v)-y.(x.v).\hspace{206
pt}\mbox{\qed}\\
\end{eqnarray*}

\vspace{0.15in}

The relation $c_{x,y}^\rho v_0=[f_x,f_y]v_0-f_{[x,y]}v_0$ (Lemma \ref{lm2.3} and the discussion after Proposition \ref{Kvzerocrit}) gives a way to calculate $c(x,y)$ explicitly.  In the notation of \S 2, $x.w_\ga=\sum_{\gamma\in\I}x_\gamma^\ga w_\gamma$ for every $x\in\fg$ and $\ga\in\I$, and

{\theorem\label{thm4.3} $\displaystyle{c(x,y)=\sum_{\ga\in\J,\gamma\in\I\setminus\J}x_\ga^\gamma y_\gamma^\ga\ \ -\sum_{\ga\in\I\setminus\J,\gamma\in\J}x_\ga^\gamma y_\gamma^\ga.}$}

\vspace{0.1in}
\noindent
\proof  Let $\kappa:\ \Vres\rightarrow {\mathbb K}$ be the linear map defined by $$\kappa (\wab w_{\ol{\gb}}^* v_0)=\left\{\begin{array}{ll}
1 & \mbox{if }\ol{\ga}=\ol{\gb}=\emptyset \\
0 & \mbox{otherwise.} \\
\end{array}
\right.
$$

\noindent
for $\wab w_{\ol{\gb}}^*\in\A^-$.  Since $c_{x,y}^\rho$ is a scalar, $c_{x,y}^\rho v_0=[f_x,f_y]v_0-f_{[x,y]}v_0=\kappa([f_x,f_y]v_0)v_0-\kappa(f_{[x,y]}v_0)v_0$.  For any $z\in\fg$, we have \begin{eqnarray*}
\kappa(f_zv_0)&=&\kappa\Big(\sum_{\ga,\gamma\in\I}z_\gamma^\ga:w_\gamma\wa^*:v_0\Big)\\
&=&\kappa\Big(\sum_{\ga\in\I}z_\ga^\ga:w_\ga\wa^*:v_0\Big)\ \ =\ \ 0,\\
\end{eqnarray*} 

\noindent
by normal ordering since $\wa\in\A^+$ or $\wa^*\in\A^+$.  Thus \begin{eqnarray*}c_{x,y}^\rho&=&\kappa\big([f_x,f_y]v_0\big)\\
&=&\kappa\Big(\displaystyle{\sum_{\ga,\gamma\in\I}y_\gamma^\ga[f_x,:w_\gamma\wa^*:]v_0}\Big)\\
&=&\kappa\Big(\displaystyle{\sum_{\ga,\gamma\in\I}y_\gamma^\ga([f_x,w_\gamma]\wa^*+w_\gamma[f_x,\wa^*])v_0}\Big)\\
&=&\kappa\Big(\displaystyle{\sum_{\ga,\gamma\in\I}x_\ga^\gamma y_\gamma^\ga(\wa\wa^*-w_\gamma w_\gamma^*)v_0}\Big).\\
\end{eqnarray*}

\noindent
Let $\chi:\ \I\rightarrow\{ 0,1\}$ with 
\begin{equation*}\label{chi}
\chi(\eta)=\left\{\begin{array}{ll}
1 & \mbox{if }\eta\in\J \\
0 & \mbox{otherwise.}
\end{array}
\right.
\end{equation*}
Note that $\wa\in\A^+$ or $\wa^*\in\A^+$ (and likewise for $w_\gamma$ and $w_\gamma^*$).  
Therefore $$\wa\wa^*v_0=
\left\{\begin{array}{ll}
0 & \mbox{if }\wa^*\in\A^+ \\
-\rho\wa^*\wa v_0+\rho v_0=\rho v_0 & \mbox{if } \wa^*\notin\A^+. \\
\end{array}
\right.$$

\noindent
Hence $(\wa\wa^*-w_\gamma w_\gamma^*)v_0=\rho(\chi(\ga)-\chi(\gamma))v_0$, so

\begin{eqnarray*} 
c(x,y)\ \ =\ \ \rho c_{x,y}^\rho&=&\rho\kappa\left(\displaystyle{\sum_{\ga,\gamma\in\I}x_\ga^\gamma y_\gamma^\ga\rho(\chi(\ga)-\chi(\gamma))v_0}\right)\\
&=&\displaystyle{\sum_{\ga,\gamma\in\I}x_\ga^\gamma y_\gamma^\ga(\chi(\ga)-\chi(\gamma))}\\
&=&\sum_{\ga\in\J,\gamma\in\I\setminus\J}x_\ga^\gamma y_\gamma^\ga -\sum_{\ga\in\I\setminus\J,\gamma\in\J}x_\ga^\gamma y_\gamma^\ga,\\
\end{eqnarray*}

\noindent
and both of these sums are finite by the finiteness condition on realizations.\qed

\section{Embeddings into Lie Algebras of Infinite \newline Matrices}

Date-Jimbo-Kashiwara-Miwa \cite{DJKM} and Kac-Peterson \cite{KP} have introduced the Lie algebra $\overline{\fA}_\infty$ consisting of infinite matrices $(a_{ij})_{i,j\in\Z}$ with $a_{ij}\in {\mathbb K}$ and only a finite number of nonzero diagonals.  The Lie algebra $\overline{\fA}_\infty$ has a central extension $\fA_\infty=\overline{\fA}_\infty\oplus {\mathbb K}\bc$ given by the 2-cocycle $\ga(E_{ij},E_{ji})=-\ga(E_{ji},E_{ij})=1$ for $i\leq 0,\ j>0$, and $\ga(E_{ij},E_{mn})=0$ otherwise (cf. \cite{KR}).

By analogy, for any realization $(\fg,W,\J,\I,\I',\rho)$, we may view the space $$\overline{\fA}:=\Big\{\displaystyle{\sum_{\ga,\gb\in\I}c_{\ga\gb}:\wa w_\gb^*:\in\Azh\ |\ c_{\ga\gb}\in {\mathbb K}}\Big\}$$  

\noindent
as a Lie algebra of $(\hbox{{\rm card}}\ \I)\times(\hbox{{\rm card}}\ \I)$ ``matrices'' $\displaystyle{\sum_{\ga,\gb\in\I}c_{\ga\gb}E_{\ga\gb}}$, where $E_{\ga\gb}$ is the matrix with $1$ in the $(\ga,\gb)$ position and $0$ elsewhere.  The Lie bracket in $\overline{\fA}$ is obtained by extending the ordinary ``matrix Lie bracket'' $[E_{\ga\gb},E_{\gamma\eta}]=\gd_{\gb,\gamma}E_{\ga\eta}-\gd_{\ga,\eta}E_{\gamma\gb}$ to $\overline{\fA}$.  Although the matrices in $\overline{\fA}$ may be of arbitrary dimension and need not have only finitely many nonzero ``diagonals,'' this bracket is well-defined, since it is simply the result of restricting the Lie bracket on $\Azh$ to $\overline{\fA}\oplus {\mathbb K}1\subset \Azh$ and then projecting onto $\overline{\fA}$.

It is clear that the space $\fA=\overline{\fA}\oplus {\mathbb K}1$, under the restriction of the Lie bracket of $\Azh$, is a central extension of $\overline{\fA}$, and the resulting 2-cocycle is simply the ``constant term'' that occurs in a given bracket---that is, $c(E_{\ga\gb},E_{\gamma\eta})=\ep\big([:\wa w_\gb^*:,:w_\gamma w_\eta^*:]\big)=\gd_{\ga,\eta}\gd_{\gb,\gamma}\big(\chi(\ga)-\chi(\gamma)\big)$, in the notation of \S 3 and \S 4. The restriction of $c$ to the Lie algebra of Corollary \ref{cor4.1} is the 2-cocycle of \S 4. If $\fg=\overline{\fA}_\infty$, $W$ is the natural representation on (doubly infinite) column vectors (with canonical basis indexed by $\I=\Z$), $\I'=\emptyset$, and $\J$ is the nonpositive integers, then we get the 2-cocycle $\ga$ and the Lie algebra $\fA_\infty$ of \cite{DJKM} and \cite{KP}.

In particular, if $V$ is faithful, the identification of the operator $:\wa w_\gb^*:$ with the matrix $E_{\ga\gb}$ and $\bc$ with $\rho 1$ gives an embedding of $\fgtilde=\fg\oplus {\mathbb K}\bc$ into $\fA$.  The following proposition gives an easy criterion for faithfulness.

{\proposition The $\fgtilde$-module $V=V(\fg,W,\J,\I,\I',\rho)$ is faithful if and only if $W$ is a faithful $\fg$-module.}

\vspace{0.1in}
\noindent
\proof  Suppose $W$ is not faithful.  Then there is a nonzero $x\in\fg$ such that $x.w=0$ for all $w\in W$.  Thus $f_x=\displaystyle{\sum_{\ga\in\I}:(x.\wa)\wa^*:=0}$, so for every $v\in V$, $x.v=f_xv=0$.  Hence $V$ is not faithful.

Conversely, suppose $W$ is faithful and $k\bc+x$ acts as zero on $V$ for some $k\in {\mathbb K}$ and $x\in\fg$.  Then for every $a\in\A$, 

\begin{eqnarray*}
0\ \ =\ \ (k\bc+x).(a v_0)&=&(k\rho+f_x)a v_0\\
&=&[k\rho+f_x,a]v_0+a(k\rho+f_x)v_0\\
&=&[f_x,a]v_0+a(k\bc+x).v_0\\
&=&[f_x,a]v_0.
\end{eqnarray*}

\noindent
Therefore by Lemma \ref{lm2.2}, $(x.\wa)v_0=0=(x.\wa^*)v_0$ for all $\ga\in\I$.  By Proposition \ref{V_res}, it now follows that $x.\wa=\displaystyle{\sum_{\gamma\in\J}x_\gamma^\ga w_\gamma}$.  Then for any $\ga,\gb\in\I$, 

\begin{eqnarray*}
0\ \ =\ \ (k\bc+x).(\wa w_\gb^*v_0)&=&[f_x,\wa w_\gb^*]v_0\\
&=&[f_x,\wa]w_\gb^*v_0+\wa[f_x,w_\gb^*]v_0\\
&=&(x.\wa) w_\gb^*v_0+\wa(x.w_\gb^*)v_0\\
&=&(x.\wa)w_\gb^*v_0\\
&=&\displaystyle{\sum_{\gamma\in\J}x_\gamma^{\ga}w_\gamma w_\gb^*v_0}\\
&=&\displaystyle{-\rho\sum_{\gamma\in\J}x_\gamma^\ga w_\gb^*w_\gamma v_0+\sum_{\gamma\in\J}x_\gamma^\ga\{ w_\gamma,w_\gb^*\rhobr v_0}\\
&=&\rho x_\gb^\ga\chi(\gb)v_0,\\
\end{eqnarray*}

\noindent
in the notation of the proof of Theorem \ref{thm4.3}.  Hence $x_\gb^\ga=0$ for all $\ga\in\I$ and $\gb\in \J$.  But $x.\wa=\displaystyle{\sum_{\gamma\in\J}x_\gamma^\ga w_\gamma}$, so this gives $x.\wa=0$ for all $\ga\in\I$.  Since $W$ is faithful, $x=0$, so $k\bc=k\bc+x$ acts as zero on $V$.  Hence $k=0$, and $V$ is faithful.\qed

\vspace{0.15in}

We can use the ``identity matrix'' $J:=\displaystyle{\sum_{\gamma\in\I}:w_\gamma w_\gamma^*:}$ to decompose $V$ into submodules:

{\proposition The $\fgtilde$-module $V$ has a $\Z$-grading $V=\bigoplus_{r\in\Z}V_r$, where $V_r$ is the $J$-eigenspace with eigenvalue $r$.  Each $V_r$ is a $\fgtilde$-submodule of $V$.  Explicitly,
$V_r=\hbox{{\rm Span}}_{\mathbb K}\{\wab\lbb\wgbs v_0\ |\ \ell(\ol{\ga})-\ell(\ol{\gb})-\ell(\ol{\gamma})=r\}.$}

\vspace{0.1in}
\noindent
\proof  For any $\gl\in W'$, $[J,\gl]=\displaystyle{\sum_{\gamma\in\I}[:w_\gamma w_\gamma^*:,\gl]}=-\rho \displaystyle{\sum_{\gamma\in\I}\{ w_\gamma,\gl\rhobr w_\gamma^*}=\displaystyle{-\sum_{\gamma\in\I}\gl(w_\gamma) w_\gamma^*}$, and for any $\ga\in\I$, $\displaystyle{-\sum_{\gamma\in\I}\gl(w_\gamma) w_\gamma^*}(\wa)=-\gl(\wa)$, so $[J,\gl]=-\gl$.  Likewise, for $w\in W$, $[J,w]=w$, and $Jv_0=0$.  Thus 

\begin{eqnarray*}
J\wab\lbb\wgbs v_0&=&[J,\wab\lbb\wgbs]v_0+\wab\lbb\wgbs Jv_0\\
&=&[J,\wab\lbb\wgbs]v_0\\
&=&\big(\ell(\ol{\ga})-\ell(\ol{\gb})-\ell(\ol{\gamma})\big)\wab\lbb\wgbs v_0.\\
\end{eqnarray*}

\noindent
To see that $V_r$ is a $\fgtilde$-submodule, it is enough to observe that $J$ commutes with $f_x$ and $\bc$ for all $x\in\fg$.  But this is trivial, since $[J,:\wa w_\gb^*:]= [J,\wa] w_\gb^*+ \wa[J, w_\gb^*]=\wa w_\gb^*-\wa w_\gb^*=0$ for all $\ga,\gb\in\I$.\qed

\vspace{0.15in}

As we will see in Example \ref{ex6.1}, the modules $V_r$ need not be irreducible, even if $W$ is irreducible.

\section{Examples}

In this section, we show that our notion of bosonic and fermionic realizations includes representations on symmetric and skew-symmetric tensors, as well as generalizes the Fock space constructions given by Feingold and Frenkel \cite{FF} for types $A_\ell$ and $A_\ell^{(1)}$ and by Gao \cite{G} for $\widehat{\mathfrak{gl}_N({\mathbb C}_q)}$, where ${\mathbb C}_q$ is the quantum torus in two variables.  Also, Theorem \ref{thm4.3} may be used to compute interesting central extensions.  We illustrate this with the Virasoro Lie algebra and a central extension of $\glA$, where $A$ is the ring of differential operators on the punctured plane ${\mathbb C}^\times$.  We anticipate that the techniques of this paper may be used to produce nontrivial representations of toroidal and other interesting Lie algebras, and we plan to investigate these and further applications in a later work.

\vspace{0.15in}

\noindent
{\example\label{ex6.1} (Symmetric and skew-symmetric tensors)}

Let $\fg$ be any Lie algebra with a module $W$.  Consider a realization $(\fg, W, \J, \I, \I')$ with $\J=\emptyset$.  Then $\A^-$ is generated by $\{ \wa\ |\ \ga\in\I\}$ and $f_xv_0=0$ for every $x\in\fg$.  The defining relations of $\A$ make it clear that $V$ is the module $S(W)$ of symmetric tensors in the tensor algebra $T(W)$ if $\rho=-1$, and is the module of skew-symmetric tensors $\bigoplus_{r\geq 1}\bigwedge^r W$, if $\rho=1$.  The central extension $\fgtilde=\fg\oplus {\mathbb K}\bc$ is obviously split, since the 2-cocycle $c$ is trivial.  Moreover, the submodules 
$$\displaystyle{V_r\cong\left\{\begin{array}{ll}
S^r(W),\mbox{\ the $r$th symmetric power of\ }W & \mbox{if\ }\rho=-1 \\
\bigwedge^r W & \mbox{if\ }\rho=1 \\
\end{array}
\right.}$$

\noindent
are seldom irreducible, even if $W$ is irreducible.

For instance, take $\fg=\mathfrak{sl}_2({\mathbb C})$, $W=S^m({\mathbb C}^2)$, $\J=\emptyset$, and $\rho=-1$, where ${\mathbb C}^2$ is the natural module and $m>1$.  Then $W$ is irreducible, but the Clebsch-Gordan rule (see \cite{FH}, for instance) gives $V_2=S^2(S^m({\mathbb C}^2))\cong\bigoplus_{0\leq n\leq\frac{m}{2}} S^{2m-4n}({\mathbb C}^2)$, so $V_2$ is \underline{not} irreducible.  Similar arguments can be made if $\rho=1$.

\vspace{0.15in}
\noindent
{\example\label{ex6.2} (Oscillator and spinor representations for $\mathfrak{gl}_N({\mathbb K})$)}

The natural representation $W={\mathbb K}^N$ of the Lie algebra $\fg=\glK$, with $\I=\{1,2,\ldots,N\}$ and $\J=\emptyset$ gives the usual oscillator and spinor representations (as described in \cite{FF}, for example) for $\rho=-1$ and $\rho=1$, respectively.  These are representations that come from the natural isomorphism $\glK\cong W\otimes W^*$, and can be thought of as the model for all bosonic and fermionic representations described in this paper (as discussed in \S 5).

\vspace{0.15in}
\noindent
{\example\label{ex6.3} (Matrices over associative algebras)}

Let $A$ be an associative algebra over the field ${\mathbb K}$, and let $M$ be an $A$-module.  The Lie algebra $\mathfrak{gl}_N(A)$ of $N\times N$ matrices over $A$ has a natural ``left-multiplication'' action on the space $M^N$ of $N\times 1$ column vectors with entries in $M$.  Specifically, let $x(a)=(x_{ij}a)_{1\leq i,j\leq N}$ where $x=(x_{ij})_{1\leq i,j\leq N}\in \mathfrak{gl}_N({\mathbb K})$ and $a\in A$, and let $v(m)=\left(\begin{array}{c}

v_1m \\
v_2m \\
\vdots \\
v_Nm
\end{array}
\right)$ where $v=\left(\begin{array}{c}

v_1 \\
v_2 \\
\vdots \\
v_N
\end{array}
\right)\in {\mathbb K}^N$ and $m\in M$.  Then $x(a).v(m)=xv(a.m)$.  Realizations of $(\mathfrak{gl}_N(A), M^N, \rho)$ give representations for a variety of interesting Lie algebras, including Examples \ref{ex6.4}, \ref{ex6.5}, and \ref{ex6.7} below.

\vspace{0.15in}
\noindent
{\example\label{ex6.4} (Affine Lie algebra $A_{N-1}^{(1)}$)}

In the notation of Example \ref{ex6.3}, let $A={\mathbb K}[t,t^{-1}]$, and let $M={\mathbb K}[t,t^{-1}]$ be the (left) regular $A$-module.  Then $W=M^N$ is a module for $\fg=\mathfrak{gl}_N(A)$.  Let $\I=\{1,2\ldots,N\}\times\Z$ with $w_{(i,n)}=e_i(n)$ and $\J=\{1,2,\ldots, N\}\times\mathbb{N}$, where $e_i(n)=e_i(t^n)$ and ${\mathbb N}$ denotes the nonnegative integers.  Then $\I'=\emptyset$, and the realization $(\fg ,W,\J,\I,\I',\rho)$ is (with the exception of a small change to the normal order)\footnote{The normal ordering used by Feingold and Frenkel is $:e_i(m)e_j(n)^*:$
$$=\frac{1}{2}e_i(m)e_j(n)^*-\frac{\rho}{2}e_j(n)^*e_i(m)=\left\{\begin{array}{ll}
-e_j(n)^*e_i(m)-\frac{\rho}{2} & \mbox{if\ }i=j\mbox{\ and\ }m=n \\
e_i(m)e_j(n)^* & \mbox{otherwise.}
\end{array}
\right.
$$  

\noindent
Their change in normal order amounts only to replacing $E_{ii}(0)$ with $E_{ii}(0)-\frac{1}{2}\bc$ in our representation.  They (and Gao \cite{G}) use the notation $a_j^*(-n)$ to denote what we call $e_j(n)^*$.} the Fock space construction given by Feingold and Frenkel for the affine Lie algebra $\fgtilde=\fg\oplus {\mathbb K}\bc$ in \cite{FF}.  If $\J=\{ 1,2,\ldots, N\}\times \Z^+$ where $\Z^+$ is the positive integers, then we recover the representation given by Gao (\cite{G}, 2.25) for his ``vertical Lie algebra'' $\mathcal{L}_v\cong\fgtilde$.

\vspace{0.15in}
\noindent
{\example\label{ex6.5} {\rm \Large(}$\widehat{\mathfrak{gl}_N(\Cq)}${\rm \Large)}}

Let $\Cq={\mathbb C}[x^\pm,y^\pm]_{n.c.}/(yx-qxy)$ be the quantum torus in two variables.\footnote{Here $q\in{\mathbb C}^\times$, and ${\mathbb C}[x^\pm,y^\pm]_{n.c.}$ is the space of Laurent polynomials in two noncommuting variables $x$ and $y$.}  The associative algebra $\Cq$ admits a representation on $M={\mathbb C}[t,t^{-1}]$ by $x.t^\ell=t^{\ell+1}$ and $y.t^\ell=q^\ell t^\ell$.  Taking $w_{(i,n)}$, $\I$, $\I'$, and $W$ as in Example \ref{ex6.4}, we let $\J=\{1,2,\ldots ,N\}\times \Z^+$.  This realization gives Gao's representation $E_{ij}(q^mp^n)\mapsto f_{ij}(m,n)=\sum_{s\in\Z}q^{ns}:e_i(s+m)e_j(s)^*:$.  In order to have the desired commutation relations between the operators $f_{ij}(m,n)$, Gao adjusts his basis of the central extension spanned by $f_{ij}(m,n)$ and $1$ by replacing the $f_{ij}(m,n)$ with $F_{ij}(m,n)=f_{ij}(m,n)+k(i,j,m,n)1$ for some constants $k(i,j,m,n)\in{\mathbb C}$.  Note that such a change of basis is included in the boundary map $\partial:\ \mathcal{C}^1\rightarrow \mathcal{C}^2$ of the Cartan-Eilenberg complex for $H^*(\mathfrak{gl}_N(\Cq);{\mathbb C})$, so does not alter the 2-cocycle or the central extension described above.\footnote{Gao's central elements $c(n)$ all act as our central element $c$ and his $c_y$ acts as $0$.}

\vspace{0.15in}
\noindent
{\example\label{ex6.6} (Virasoro Lie algebra)}

Let $\fg$ be the Witt algebra $\hbox{{\rm Span}}_{\mathbb K}\{ L_m\ |\ m\in \Z\}$ with commutation relations $[L_m,L_n]=(m-n)L_{m+n}$ for all $m,n\in\Z$.  The Witt algebra has a natural (and faithful) representation as the derivations $L_m=-t^{m+1}\frac{d}{dt}$ on the space of Laurent polynomials, $W={\mathbb K}[t,t^{-1}]$.  Taking $\I=\Z$, $\I'=\emptyset$, and $\J=\Z^+$, we obtain the 2-cocycle 
$$c(L_m,L_n)=\sum_{i\geq 0,j<0}(L_m)_i^j(L_n)_j^i-\sum_{i< 0,j\geq 0}(L_m)_i^j(L_n)_j^i.$$  

\noindent
Clearly $L_m.t^j=-jt^{m+j}$, so $(L_m)_i^j=-j\gd_{i,m+j}$ and 

\begin{eqnarray*}
c(L_m,L_n)&=&\sum_{i\geq 0,j<0}ij\gd_{i,m+j}\gd_{j,n+i}-\sum_{i< 0,j\geq 0}ij\gd_{i,m+j}\gd_{j,n+i}\\
&=&\sum_{i=0}^{m-1}\gd_{m,-n}i(i-m)-\sum_{i=m}^{-1}\gd_{m,-n}i(i-m)\\
&=&\gd_{m,-n}\left(\frac{m-m^3}{6}\right).\\
\end{eqnarray*}

\noindent
Adjusting our 2-cocycle $c$ by a constant factor of $-\frac{1}{2}$ gives the usual normalization of the central extension $\fgtilde=\fg\oplus {\mathbb K}\bc$, called the {\it Virasoro Lie algebra}.  It is well-known that this central extension is universal (cf. \cite{BM}, \cite{KR}, or \cite{MP}, for instance).

\vspace{0.15in}
\noindent
{\example\label{ex6.7} {\rm \Large(}$\widehat{\mathfrak{gl}_N}(\mathcal{W}_{1+\infty})${\rm \Large)}}

Let ${\mathbb K}[p,q^{\pm 1}]_{n.c.}$ be the (unital) associative algebra of noncommuting polynomials in two variables $p$ and $q$, localized at the ideal $(q)$.  Let $A={\mathbb K}[p,q^{\pm 1}]_{n.c.}/([pq]-1)$.  The algebra $A$ is the first Weyl algebra $A_1$ localized on the multiplicative set $\{q^m\ |\ m\in\mathbb{N}\}$, and thus has a natural representation on ${\mathbb K}[t,t^{-1}]$ given by $
p.t^\ell=\ell t^{\ell-1},\  
q.t^\ell=t^{\ell+1},\mbox{ and }  
q^{-1}.t^\ell=t^{\ell-1}.$  If ${\mathbb K}$ is algebraically closed, $A$ is also the ring of differential operators on the ``circle'' $\{ (x,y)\in {\mathbb K}^2\ |\ x^2+y^2=1\}$ or the scheme ${\mathbb K}\setminus\{ 0\}$ (cf. \cite{Bj}, \cite{Ha}).\footnote{which is the punctured plane in the case ${\mathbb K}={\mathbb C}$.}  Viewed as a Lie algebra in the usual way, the algebra $A$ arises as a limit of objects called $W_n$-algebras (see \cite{Z}, for instance), and is usually denoted $\mathcal{W}_{1+\infty}$. 

In the notation of Example \ref{ex6.3}, $W=({\mathbb K}[t,t^{-1}])^N$ is a representation for the Lie algebra $\fg=\mathfrak{gl}_N(A)$, and we obtain a realization by taking $w_{(i,n)}=e_i(t^n)\in W$, $\I=\{ 1,2,\ldots, N\}\times\Z$, $\I'=\emptyset$, and $\J=\{ 1,2,\ldots, N\}\times {\mathbb N}$.  Let $x(k,\ell)=x(q^kp^\ell)$ where $x\in\mathfrak{gl}_N({\mathbb K})$, $k\in\Z$, and $\ell\in\mathbb{N}$.  Since the module action is $E_{ij}(k,\ell).w_{(m,n)}=\gd_{j,m}(n)_\ell w_{(i,n+k-\ell)}$ where $(n)_\ell=n(n-1)\cdots (n-(\ell-1))$, we have $E_{ij}(k,\ell)_{(r,s)}^{(m,n)}=\gd_{i,r}\gd_{j,m}\gd_{n+k-\ell,s}(n)_\ell$.  Thus the Fock space representation is $$E_{ij}(k,\ell)\mapsto \sum_{n\in\Z}(n)_{\ell}:w_{(i,n+k-\ell)}w_{(j,n)}^*:,$$ and by Theorem \ref{thm4.3}, the resulting 2-cocycle is

\bigskip
\noindent
$\displaystyle{c(E_{ij}(k,\ell),E_{mn}(r,s))}$
\begin{eqnarray*}
&=&\sum_{(a,b)\in\J,(u,v)\in\I\setminus\J}E_{ij}(k,\ell)_{(a,b)}^{(u,v)}E_{mn}(r,s)_{(u,v)}^{(a,b)}\\
&&\hspace{0.5in}-\sum_{(a,b)\in\I\setminus\J,(u,v)\in\J}E_{ij}(k,\ell)_{(a,b)}^{(u,v)}E_{mn}(r,s)_{(u,v)}^{(a,b)} \\
&=& \gd_{i,n}\gd_{j,m}\Big(\sum_{b\geq 0,v<0}\gd_{v+k-\ell,b}\gd_{b+r-s,v}(v)_\ell(b)_s\\
&&\hspace{0.5in}-\sum_{b<0,v\geq 0}\gd_{v+k-\ell,b}\gd_{b+r-s,v}(v)_\ell(b)_s\Big)\\
&=&\gd_{k+r,\ell+s}\gd_{j,m}\gd_{i,n}\Big(\sum_{b=0}^{s-r-1}(b+r-s)_\ell(b)_s-\sum_{b=s-r}^{-1}(b+r-s)_\ell(b)_s\Big)\\
\end{eqnarray*}

\noindent
with the convention that empty sums are $0$.

By the lemma and identities in the Appendix, if $s\geq r$, then 
\begin{eqnarray*}
\sum_{b=0}^{s-r-1}(b+r-s)_\ell(b)_s&=&(-1)^\ell\ell !s!\binom{s-r+\ell}{-r-1} \\
&=&(-1)^\ell\ell!s!\binom{s-r+\ell}{s+\ell+1}\\
&=&(-1)^{s+1}\ell!s!\binom{r}{s+\ell+1}. \\
\end{eqnarray*}

\noindent
Likewise, if $s<r$, then $r>0$, so 

\begin{eqnarray*}
\sum_{b=s-r}^{-1}(b+r-s)_\ell(b)_s&=&(-1)^s\ell!s!\binom{r}{r-s-\ell-1}\\
&=&(-1)^s\ell!s!\binom{r}{s+\ell+1}.\\
\end{eqnarray*}

\noindent
Hence $c(E_{ij}(k,\ell),E_{mn}(r,s))$
$$\begin{array}{l}
=\ \gd_{k+r,\ell+s}\gd_{j,m}\gd_{i,n}(-1)^{s+1}\ell!s!\binom{r}{s+\ell+1}\\
=\ \hbox{{\rm tr}}(E_{ij}E_{mn})\gd_{k+r,\ell+s}(-1)^{s+1}\ell!s!\binom{r}{s+\ell+1},\\
\end{array}$$

\noindent
so by the linearity of the trace, $$c(x(k,\ell),y(r,s))=\hbox{{\rm tr}}(xy)\gd_{k+r,\ell+s}(-1)^{s+1}\ell!s!\binom{r}{\ell+s+1}$$ 

\noindent
for any $x,y\in\mathfrak{gl}_N({\mathbb K})$.

The 2-cocycle $c$ was first computed by Kac-Peterson \cite{KP}, and the corresponding central extension $\widehat{\mathfrak{gl}_N(A)}$ is universal (\cite{Fei}, \cite{La}, \cite{Li}).  With a slight abuse of notation, we denote the central extension as $\widehat{\mathfrak{gl}_N}(\mathcal{W}_{1+\infty})$.  By the 2-cocycle formula obtained above, $\widehat{\mathfrak{gl}_N}(\mathcal{W}_{1+\infty})$ contains full copies of the affine and Virasoro Lie algebras:
$$\widehat{\mathfrak{gl}_N({\mathbb K})}\cong \hbox{{\rm Span}}_{\mathbb K}\{ x(k,0)\ |\ x\in\mathfrak{gl}_N({\mathbb K}),\ k\in\Z\}\oplus {\mathbb K}\bc$$
$$Vir\cong \hbox{{\rm Span}}_{\mathbb K}\{ I(k,1)\ |\ k\in\Z\}\oplus {\mathbb K}\bc$$
\noindent
where $I\in\mathfrak{gl}_N({\mathbb K})$ is the identity matrix.

\section{Appendix: A Combinatorial Lemma}

In the statements below, we assume that all variables represent integers, and we follow the standard conventions that empty sums are $0$, empty products are $1$, and $\displaystyle{\binom{a}{b}=\frac{(a)_b}{b!}}$ where
$$(a)_b=\left\{\begin{array}{ll}

a(a-1)\cdots (a-(b-1)) & \mbox{ if }b\geq 0\\
0 & b<0.
\end{array}
\right.$$

\noindent
We will use the following well-known identities without explicit mention:

\begin{itemize}
\item[{\rm (I1)}] $\displaystyle{\binom{a}{b}=(-1)^b\binom{b-a-1}{b}}$
\item[{\rm (I2)}] $\displaystyle{\sum_{n=0}^c\binom{a}{n}\binom{b}{c-n}=\binom{a+b}{c}}$
\item[{\rm (I3)}] $\displaystyle{\binom{d}{b}=\binom{d}{d-b}}$
\end{itemize}

\noindent
for any $a,b,c\in\Z$ and $d\in\mathbb{N}$.  (See \cite{Fe}, for instance.)

{\lemma Let $m,\ell,s\geq 0$.  Then 

\begin{itemize}
\item[{\rm (i)}]$\displaystyle{\sum_{b=0}^{m-1}(b-m)_\ell(b)_s=(-1)^\ell\ell!s!\binom{m+\ell}{m-s-1}}$
\item[{\rm (ii)}]$\displaystyle{\sum_{b=-m}^{-1}(b+m)_\ell(b)_s=(-1)^s\ell!s!\binom{m+s}{m-\ell-1}}$.
\end{itemize}
}

\vspace{0.1in}
\noindent
\proof  Note that (ii) will follow from (i) by interchanging $\ell$ and $s$ and making the change of variables $b\mapsto b+m$ in (i).  Thus we may restrict our attention to the first identity:

\begin{eqnarray*}
\sum_{b=0}^{m-1}(b-m)_\ell(b)_s&=&\sum_{b=0}^{m-1}(-1)^\ell(m-b+\ell-1)_\ell(b)_s\\
&=&\sum_{b=s}^{m-1}(-1)^\ell\ell!s!\binom{m-b+\ell-1}{\ell}\binom{b}{s}\\
&=&(-1)^\ell\ell!s!\sum_{b=s}^{m-1}\binom{m-b+\ell-1}{m-b-1}\binom{b}{b-s}\\
&=&(-1)^{\ell+m-s-1}\ell!s!\sum_{b=s}^{m-1}\binom{-\ell-1}{m-b-1}\binom{-s-1}{b-s}.\\
\end{eqnarray*}

\noindent
Changing the index of summation from $b$ to $b-s$ now gives

\begin{eqnarray*}
\sum_{b=0}^{m-1}(b-m)_\ell(b)_s&=&(-1)^{\ell+m-s-1}\ell!s!\sum_{b=0}^{m-1-s}\binom{-\ell-1}{m-b-s-1}\binom{-s-1}{b}\\
&=&(-1)^{\ell+m-s-1}\ell!s!\binom{-\ell-s-2}{m-s-1}\\
&=&(-1)^\ell\ell!s!\binom{m+\ell}{m-s-1}.\\
\end{eqnarray*}
\qed

\end{document}